\documentclass{article}

\usepackage{a4wide}
\usepackage{amsmath,amssymb,amsthm}
\usepackage{txfonts}
\usepackage{mathtools}
\usepackage{xcolor}
\usepackage{hyperref}
\usepackage{mdwlist}

\newtheorem{remark}{Remark}
\newtheorem{lemma}{Lemma}

\newtheorem{theorem}{Theorem}

\renewcommand{\eqref}[1]{\hyperref[#1]{(\ref{#1})}}

\makeatletter
\newcommand*\rel@kern[1]{\kern#1\dimexpr\macc@kerna}
\newcommand*\widebar[1]{%
  \begingroup
  \def\mathaccent##1##2{%
    \rel@kern{0.8}%
    \overline{\rel@kern{-0.8}\macc@nucleus\rel@kern{0.2}}%
    \rel@kern{-0.2}%
  }%
  \macc@depth\@ne
  \let\math@bgroup\@empty \let\math@egroup\macc@set@skewchar
  \mathsurround\z@ \frozen@everymath{\mathgroup\macc@group\relax}%
  \macc@set@skewchar\relax
  \let\mathaccentV\macc@nested@a
  \macc@nested@a\relax111{#1}%
  \endgroup
}
\makeatother

\newcommand{\eps}{\varepsilon}
\newcommand{\RR}{\mathbb{R}}
\newcommand{\NN}{\mathbb{N}}
\renewcommand{\phi}{\varphi}
\renewcommand{\rho}{\varrho}
\renewcommand{\theta}{\vartheta}

\newcommand{\m}[1]{\mathcal{#1}}

\DeclareMathOperator{\ess}{ess}

\newcommand{\E}{\mathrm{e}}
\newcommand{\D}{\mathrm{  d}}

\newcommand{\abs}[1]{\left|#1\right|}
\newcommand{\norm}[1]{\left\|#1\right\|}
\newcommand{\bignorm}[1]{\big\|#1\big\|}

\hypersetup{
  pdftitle    = {Maximum-norm \emph{a posteriori} error bounds for an extrapolated upwind scheme
  applied to a singularly perturbed convection-diffusion problem},
  pdfsubject  = {Numerical Analysis},
  pdfauthor   = {Torsten Lin\ss, Goran Radojev},
  pdfkeywords = {maximum-norm a posteriori error estimates, extrapolation,
  upwind sheme, singular perturbation, convection-diffusion, adaptivity},
  pdfcreator  = {pdflatex},
  pdfproducer = {LaTeX with hyperref},
    colorlinks,
    linkcolor={red!90!black},
    citecolor={blue!90!black},
    urlcolor={green!90!black}
}

\title{Maximum-norm \emph{a posteriori} error bounds for an extrapolated upwind scheme
  applied to a singularly perturbed convection-diffusion problem}

\author{Torsten Lin\ss\thanks{Fakult\"at f\"ur Mathematik und Informatik,
        FernUniversit\"at in Hagen,
        Universit\"atsstra{\ss}e 11,
        58095 Hagen,
        Germany,
        \texttt{torsten.linss@fernuni-hagen.de}}
   \and Goran Radojev\thanks{Department of Mathematics and Computer Science, Faculty of Sciences,
        University of Novi Sad, Trg Dositeja Obradovi\'ca~4, 21000 Novi Sad,
        Serbia,
        \texttt{goran.radojev@dmi.uns.ac.rs}.
        GR acknowledges financial support from the Faculty of Mathematics and Computer Science
        at FernUniversität in Hagen through a visiting grant.}}

\begin{document}

\maketitle

\begin{abstract}
  Richardson extrapolation is applied to a simple first-order upwind difference
  scheme for the approximation of solutions of singularly perturbed
  convection-diffusion problems in one dimension.
  Robust \emph{a posteriori} error bounds are derived for the proposed method
  on arbitrary meshes.
  It is shown that the resulting error estimator can be used to stear an
  adaptive mesh algorithm that generates meshes resolving layers and
  singularities.
  Numerical results are presented that illustrate the theoretical findings.


\end{abstract}


\section{Introduction}

Extrapolation techniques are a very powerful and frequently used tool to
enhance the approximation accuracy of an underlying numerical scheme, see
\cite{MR0307435} for a historical survey.
The first application to the numerical treatment of differential equations
is due to Richardson~\cite{Richardson1910}.

In recent years there has been a number of publications dedicated to extrapolation
methods applied to singularly perturbed problem.
Some of the early ones dedicated to the convergence analysis of these methods in
the \emph{maximum norm} on \emph{a priori} chosen layer-adapted meshes
are~\cite{MR1967580, MR2046176, MR2035009}.
The situation is different when it comes to maximum-norm \emph{a posteriori} error
estimation for singularly perturbed problems.
Pioneering work was conducted by Kopteva~\cite{MR1860270,MR1870850} for first-order
upwind scheme.
This technique was extended in~\cite{MR2644301} to a second-order defect correction
method.
The aim of the present publication is its adaptation to Richardson extrapolation
for a singularly perturbed convection-diffusion problem.

Let us consider the convection-diffusion problem in conservative form
\begin{equation}\label{problem}
 \m{L}u \coloneqq -\varepsilon u'' - (bu)' + c u   = f\,,
             \ \ \text{in} \ \ (0,1), \ \ u(0)  = \gamma_0, \ u(1)=\gamma_1,
\end{equation}
where $\eps$ is a small positive parmeter and \mbox{$b\ge\beta>0$} and \mbox{$c\ge 0$}
on \mbox{$[0,1]$}.
It provides an excellent paradigm for the study of numerical techniques in
computational fluid dynamics and for the treatment of problems with boundary
layers, i.e., regions where the solution and its derivatives change
rapidly~\cite{MR2454024}.
The solution to~\eqref{problem} exhibits an exponential boundary layer of width
$\mathcal{O}(\eps\ln{1/\varepsilon})$ near the point $x=0$. 

There is a vaste literature studying numerical methods for solving this problem.
The reader is referred to~\cite{MR2454024} for a survey.
In most of these studies, the focus is on \emph{a priori} error estimation.
Only a few papers deal with \emph{a posteriori} error estimation.

The main result of the present paper, \autoref{theo:main} is an \emph{a posteriori}
error bound for Richardson extrapolation applied to a first-order upwind method when
applied to~\eqref{problem}.
This result is robust in the sense that all dependencies on the perturbation
parameter~$\eps$ are carefully tracked.
Based on these bounds we device an adaptive algorithm that is able to condense
and adapt the mesh in regions where and as necessary.

The paper is organized as follows.
In \autoref{sect:method} we describe the extrapolation procedure and the upwind
scheme upon which it is based.
The main result -- a robust a posteriori error estimate -- is derived in
\autoref{sect:analysis}.
The final section documents some results of our numerical experiments.
We test the error estimator when applied to two standard \emph{a priori} adapted
meshes.
Moreover, we adapt and test de Boor's mesh movement algorithm to automatically
adapt the mesh without a priori knowledge of the location of layers.

\section{The upwind difference scheme and extrapolation}
\label{sect:method}

On $[0,1]$ consider the mesh
\begin{gather*}
  \omega \colon 0=x_0<x_1<x_2<\cdots<x_N=1,
\end{gather*}
with $N\in\NN$ being the number of mesh intervals,
$I_i\coloneqq (x_{i-1},x_i)$, $h_i\coloneqq x_i-x_{i-1}$ and
$x_{i-1/2} \coloneqq \bigl(x_i+x_{i-1}\bigr)/2$.
For any function $\phi\in C[0,1]$ we set $\phi_i\coloneqq \phi(x_i)$
and $\phi_{i-1/2}\coloneqq \phi(x_{i-1/2})$.

Our discretisation is based on a simple first-order upwind difference scheme
on the mesh $\omega$: Find $V\in\RR^{N+1}$ such that $V_0 = \gamma_0$,
$V_N = \gamma_1$ and
\begin{subequations}\label{extra}
\begin{alignat}{2}\label{one-step}
  - \frac{\eps}{h_{i+1}}
    \left(\frac{V_{i+1}-V_i}{h_{i+1}}-\frac{V_i-V_{i-1}}{h_i}
    \right)
  - \frac{(b V)_{i+1}-(bV)_i}{h_{i+1}} + c_i V_i & = f_i, &\ \ & i=1,\dots, N-1.
  \intertext{A second approximation $W\in\RR^{2N+1}$ is computed by applying
     the upwind scheme on a mesh obtained by simply bisecting $\omega$.
     It satisfies $W_0 = \gamma_0$, $W_N = \gamma_1$ and}
  \label{two-step-1}
  - \frac{\eps}{h_{i+1}/2}
    \left(\frac{W_{i+1}-W_{i+1/2}}{h_{i+1}/2}-\frac{W_{i+1/2}-W_{i}}{h_{i+1}/2}
    \right)
  - \frac{(b W)_{i+1}-(bW)_{i+1/2}}{h_{i+1}/2} + c_{i+1/2} W_{i+1/2} & = f_{i+1/2}\,, &\ \ & i=0,\dots, N-1,\\
  \label{two-step-2}
  - \frac{\eps}{h_{i+1}/2}
    \left(\frac{W_{i+1/2}-W_{i}}{h_{i+1}/2}-\frac{W_i-W_{i-1/2}}{h_i/2}
    \right)
  - \frac{(b W)_{i+1/2}-(bW)_i}{h_{i+1}/2} + c_i W_i & = f_i\,, &\ \ & i=1,\dots, N-1.
\end{alignat}
These two are combined by means of the Richardson extrapolation procedure
to give our final extrapolated approximation:
\begin{gather}\label{comb}
  U_i = 2W_i - V_i\,, \ \ i=0,\dots, N.
\end{gather}
\end{subequations}

\section{A posteriori error analysis}
\label{sect:analysis}

\subsection{Stability of the differential operator}

Throughout, we shall assume, that the coefficients in~\eqref{problem} satisfy
\begin{gather}\label{ass-data}
  c \ge 0 \quad \text{and}\quad c-b'\ge 0 \quad \text{on} \quad [0,1]\,.
\end{gather}
\begin{remark}
 The results hold also without the restrictions imposed by~\eqref{ass-data},
 see~\cite{0999.65079}, however the arguments become more complicated and the
 constants in the estimates will differ.

 Furthermore, note that when $\beta>0$ then~\eqref{ass-data} can always
 be ensured for $\eps$ smaller than some threshold value $\eps_0$ by a simple
 transformation $u(x)=\tilde{u}(x)\E^{\chi x}$ with $\chi$ chosen appropriately.
\end{remark}

For any $v\in C^1[0,1]$ define the differential-integral operator $\m{A}$ by
\begin{gather}\label{def-A}
  \bigl(\m{A}v\bigr)(x)
    \coloneqq \eps v'(x)+ (bv)(x) + \int_x^1 \bigl(cv\bigr)(s) \D s\,,
    \ \ x\in[0,1].
\end{gather}
Furthermore, set
\begin{gather*}
 \m{F}(x) \coloneqq \int_x^1 f(s) \D s\,, \ \ x\in[0,1].
\end{gather*}
Clearly,
\begin{gather*}
  \m{L}v = - \bigl(\m{A}v\bigr)' \quad \text{and} \quad
  f = \m{F}' \quad \text{on} \ \ [0,1].
\end{gather*}
Furthermore, for the solution $u$ of~\eqref{problem} there holds
\begin{gather}\label{int-Lu-f}
  \m{A}u - \m F \equiv \alpha \quad \text{on} \ \ [0,1]\,,
   \quad \text{with} \ \ \alpha = \eps u'(1) + b(1)u(1).
\end{gather}

For our further investigations, let us introduce the supremum and the
$L_1$ norms
\begin{gather*}
  \norm{v}_\infty \coloneqq \ess \sup_{s\in[0,1]} \abs{v(s)}
 \qquad\text{and}\qquad
  \norm{v}_1 \coloneqq \int_0^1 \abs{v(s)}\D s.
\end{gather*}

We have the following stability results.
\begin{lemma}\label{lem:stab-cont}
  Suppose~\eqref{ass-data} holds true.
  Then the operator $\m{A}$ satisfies the following stability properties.
  \begin{enumerate}
    \item[(a)]
      \begin{gather*}
        \norm{v}_\infty \le \frac{2}{\beta} \min_{\kappa \in \RR} \norm{\m{A} v + \kappa}_\infty
             \ \ \  \text{for all} \ \ v\in W_0^{1,\infty}(0,1).
      \end{gather*}
    \item[(b)]
      For any function $w\in W_0^{1,\infty}(0,1)$ with
      \begin{gather*}
        \bigl(\m{A} w\bigr)(x)= \Xi_{i-1/2}(x-x_{i-1/2})
          \ \ \text{for} \ \ x\in(x_{i-1},x_i), \ \ i=1,\dots,N,
      \end{gather*}
      there holds
      \begin{gather*}
       \norm{w}_\infty \le
         C^* \max_{i=1,\dots,N}
              \left\{ \abs{\Xi_{i-1/2}}
                             \min\left[\frac{h_i}{\norm{b}_\infty},
                                       \frac{h_i^2}{4\eps}
                                 \right]
              \right\},
      \end{gather*}
      where
      \begin{gather*}
         C^*= \frac{2\norm{b}_\infty + \norm{c}_\infty + \beta}{2\beta}\,.
      \end{gather*}
  \end{enumerate}
\end{lemma}
\begin{proof}
  For (a) see~\cite{MR1937332}, while (b) was given in~\cite{MR2644301}.
\end{proof}

\subsection{Error estimation}

Using \autoref{lem:stab-cont}, we can now derive our a posterior error bound
for the extrapolation methods~\eqref{extra}.
To this end we need to extend \mbox{$U\colon \omega \to \RR$} to a function
defined on all of $[0,1]$.
We do so using linear interpolation.
For any function \mbox{$\phi\colon \omega\to \RR$}, we introduce its linear
interpolant by
%
\begin{gather*}
  \bar{\phi}(x) = \frac{\phi_{i-1}+\phi_i}{2}
                    + \bigl(x-x_{i-1/2}\bigr)\phi_{x,i},
      \ \ \ x\in[x_{i-1},x_i], \ \ i=1,\dots,N,
  \intertext{where}
    \phi_{x,k} \coloneqq \frac{\phi_k - \phi_{k-1}}{h_k}\,, \ \ k=1,\dots,N.
\end{gather*}
For functions $\phi\in C[0,1]$ we define another piecewise linear,
but (possibly) discontinuous interpolant:
\begin{gather*}
  \hat{\phi}(x) \coloneqq \phi_{i-1/2} + \bigl(x-x_{i-1/2}\bigr) \phi_{x,i}\,,
    \ \ \ x\in (x_{i-1},x_i), \ \ i=1,\dots,N.
\end{gather*}

Then, by~\eqref{int-Lu-f} and~\eqref{comb}
\begin{gather}\label{A(u-U)}
  \left(\m{A}\bigl(u-\bar{U}\bigr)\right)(x)
     = \m{F}(x) + \alpha - 2 \bigl(\m{A}\bar{W}\bigr)(x)
                         + \bigl(\m{A}\bar{V}\bigr)(x).
\end{gather}
Next, we have to derive expressions for $\m{A}\bar{V}$ and $\m{A}\bar{W}$.

Multiplying~\eqref{one-step} by $h_{i+1}$ and summing for
\mbox{$i=k,\dots,N-1$}, we get
\begin{gather}
   \eps \frac{V_k-V_{k-1}}{h_k} + (bV)_k
        + \sum_{i=k}^{N-1} h_{i+1} \left(cV -f\right)_i
    \equiv \alpha^V =
       \eps \frac{V_N-V_{N-1}}{h_N} + (bV)_N, \ \ k=1,\dots,N-1.
\end{gather}
This representation together with the definition of $\m{A}$, gives
\begin{align}\label{AbarV}
  \bigl(\m{A}\bar{V}\bigr)(x)
     & = \eps \frac{V_k-V_{k-1}}{h_k} + b(x)\bar{V}(x)
          + \int_x^1 c(s)\bar{V}(s) \D s \notag \\
     & = \bigl(b\bar{V}\bigr)(x) - b_k V_k + \int_x^1 c(s)\bar{V}(s) \D s
          - \sum_{i=k}^{N-1} h_{i+1} \left(cV-f\right)_i + \alpha^V \\
     & \hspace{20em} \text{for} \ \ x\in I_k, \ \ k=1,\dots,N-1. \notag
\end{align}
Next, multiply~\eqref{two-step-1} and~\eqref{two-step-2} by $h_{i+1}/2$ and
sum both equations for $i=k,\dots,N-1$ to verify that
\begin{multline*}
  \eps \frac{W_k-W_{k-1/2}}{h_k/2} + (bW)_k
       + \sum_{i=k}^{N-1} \frac{h_{i+1}}{2}
         \left[\left(cW-f\right)_{i+1/2} + \left(cW-f\right)_i\right] \\
   \equiv \alpha^W =
      \eps \frac{W_N-W_{N-1}}{h_N} + (bW)_N, \ \ k=1,\dots,N-1.
\end{multline*}
This yields
\begin{align*}
  2\bigl(\m{A}\bar{W}\bigr)(x)
     & = - \eps \frac{W_k - 2 W_{k-1/2} + W_{k-1}}{h_k/2}
         + 2 b(x) \bar{W}(x) - 2 b_k W_k + 2 \int_x^1 c(s)\bar{W}(s) \D s \\
     & \qquad\quad
       - \sum_{i=k}^{N-1} h_{i+1}
         \left[\left(cW-f\right)_{i+1/2} + \left(cW-f\right)_i\right] + 2\alpha^W
                \ \ \ \text{for} \ \ x\in I_k, \ \ k=1,\dots,N-1.
\end{align*}
Use~\eqref{two-step-1} with $i=k-1$ to eliminate
$\bigl(W_k - 2 W_{k-1/2} + W_{k-1}\bigr)/h_k$.
We get
\begin{gather*}
  \begin{split}
  2\bigl(\m{A}\bar{W}\bigr)(x)
     & = 2 \bigl(b\bar{W}\bigr)(x) - (bW)_k - (bW)_{k-1/2}
                + 2 \int_x^1 c(s)\bar{W}(s) \D s
               - \frac{h_k}{2} \left(cW-f\right)_{k-1/2} \\
     & \qquad\quad
       - \sum_{i=k}^{N-1} h_{i+1}
         \left[\left(cW-f\right)_{i+1/2} + \left(cW-f\right)_i\right] + 2\alpha^W
                \ \ \ \text{for} \ \ x\in I_k, \ \ k=1,\dots,N-1.
  \end{split}
\end{gather*}
Substituting~\eqref{AbarV} and the last equation into~\eqref{A(u-U)},
we arrive at
%
%
\begin{align}
  & \left(\m{A}\bigl(u-\bar{U}\bigr)\right)(x) \notag \\
  & \qquad = \int_x^1 \bigl(f-c\bar{U}\bigr)(s) \D s
       + \sum_{i=k}^{N-1} h_{i+1}
         \left[\left(cW-f\right)_{i+1/2} + \left(cW-f\right)_i
                 - \left(cV-f\right)_i \right]
               + \frac{h_k}{2} \left(cW-f\right)_{k-1/2} \notag \\
  & \label{Au-U-expanded}
    \qquad\qquad
       - \bigl(b\bar{U}\bigr)(x) + (bW)_k + (bW)_{k-1/2} - b_k V_k
       + \alpha - 2 \alpha^W + \alpha^V
                \ \ \ \text{for} \ \ x\in I_k, \ \ k=1,\dots,N-1.
\end{align}

The terms involving $b$ can be rewritten as follows.
First, a direct calculation shows that
\begin{gather*}
  \bigl(b\bar{U}\bigr)(x)
    = \bigl(b\bar{U} - \widebar{bU}\bigr)(x)
                     + \widebar{bU}(x)
    = \bigl(b\bar{U} - \widebar{bU}\bigr)(x)
                     + \widebar{bU}_{k-1/2} + \bigl(x-x_{k-1/2}\bigr) (bU)_{x,k}\,,
\end{gather*}
which yields
\begin{gather}\label{b-expanded}
\begin{split}
  & - \bigl(b\bar{U}\bigr)(x) + (bW)_k + (bW)_{k-1/2} - b_k V_k \\
  & \qquad\quad
    = \underbrace{\bigl(\widebar{bU}-b\bar{U}\bigr)(x)}_{\displaystyle \eqqcolon B(x)}
        - \bigl(x-x_{k-1/2}\bigr)(bU)_{x,k}
         + \underbrace{(bW)_{k-1/2} - (bW)_{k-1} - \frac{(bV)_k - (bV)_{k-1}}{2}}_{
                   \displaystyle \eqqcolon \Delta_k}
\end{split}
\end{gather}

For the terms in~\eqref{Au-U-expanded} that involve integrals and sums, we
proceed as follows.
Let $\psi\coloneqq f-c\bar{U}$.
Then
\begin{align}\notag
  \int_x^1 \psi(s) \D s
    & =
     \int_x^1 \bigl(\psi-\hat{\psi}\bigr)(s) \D s
       + \bigl(x_k-x\bigr) \psi_{k-1/2}
       + \int_x^{x_k} \bigl(s-x_{k-1/2}\bigr) \D s \ \psi_{x,k}
       + \sum_{i=k}^{N-1} h_{i+1} \psi_{i+1/2} \\
    & =
     \underbrace{
     \int_x^1 \bigl(\psi-\hat{\psi}\bigr)(s) \D s
       + \frac{(x_k-x)(x-x_{k-1})}{2} \ \psi_{x,k}}_{\displaystyle \eqqcolon \Psi(x)}
       + \bigl(x_{k-1/2}-x\bigr) \psi_{k-1/2}
       + \sum_{i=k}^{N-1} h_{i+1} \psi_{i+1/2}
       + \frac{h_k}{2} \psi_{k-1/2}
  \label{intpsi}
\end{align}
The last two term can be combined with the sums in~\eqref{Au-U-expanded} to
give
\begin{align} \notag
  & \sum_{i=k}^{N-1} h_{i+1} \psi_{i+1/2}
       + \frac{h_k}{2} \psi_{k-1/2}
   + \sum_{i=k}^{N-1} h_{i+1}
         \left[\left(cW-f\right)_{i+1/2} + \left(cW-f\right)_i\right]
               + \frac{h_k}{2} \left(cW-f\right)_{k-1/2}
          - \sum_{i=k}^{N-1} h_{i+1} \left(cV-f\right)_i \\
  & \quad = \sum_{i=k}^{N-1} h_{i+1} \left\{
               \bigl(c(W-\bar{U})\bigr)_{i+1/2} + \bigl(c(W-V)\bigr)_i
               \right\}
               + \frac{h_k}{2} \bigl(c(W-\bar{U})\bigr)_{k-1/2} \eqqcolon \Gamma_k
  \label{sums}
\end{align}
Then, subsituting~\eqref{b-expanded}, \eqref{intpsi} and~\eqref{sums} into
\eqref{Au-U-expanded}, we get
\begin{align*}
  \left(\m{A}\bigl(u-\bar{U}\bigr)\right)(x)
     = \Psi(x) + B(x) + \left(\Gamma+\Delta\right)_k
          + \bigl(x_{k-1/2}-x\bigr)\left\{(f-c\bar{U})_{k-1/2} + (bU)_{x,k}\right\}
       + \alpha - 2 \alpha^W + \alpha^V\,.
\end{align*}
Finally, application of \autoref{lem:stab-cont} establishes our main result.
\begin{theorem}\label{theo:main}
  Suppose~\eqref{ass-data} holds true.
  Set $\psi\coloneqq f-c\bar{U}$.
  Then the error of the extrapolated upwind scheme~\eqref{extra} satisfies
  \begin{align*}
    \norm{u-\bar{U}}_\infty
      & \le \frac{2}{\beta}
              \left\{\norm{\psi-\hat{\psi}}_1
                       + \max_{k=1,\dots,N} \frac{h_k^2}{8} \abs{\psi_{x,k}}
                       + \bignorm{\widebar{bU} - b\bar{U}}_\infty
                       + \max_{k=1,\dots,N} \abs{\Gamma_k + \Delta_k} \right\} \\
      & \qquad + C^* \max_{k=1,\dots,N}
                \left\{ \abs{(f-c\bar{U})_{k-1/2} + (bU)_{x,k}}
                               \min\left[\frac{h_k}{\norm{b}_\infty},
                                         \frac{h_k^2}{4\eps}
                                   \right]
                \right\} \eqqcolon \eta_{N}
  \end{align*}
\end{theorem}

\begin{remark}
  The error bound of~\autoref{theo:main} contains terms, namely
  \mbox{$\norm{\psi-\hat{\psi}}_1$} and
  \mbox{$\bignorm{\widebar{bU} - b\bar{U}}_\infty$} that in general have
  to be approximated.
  For example, Simpson's rule can be applied on each interval $I_k$ to obtain
  \begin{gather*}
    \norm{\psi-\hat{\psi}}_1
       \approx \sum_{k=1}^N h_k\frac{\abs{\psi_k - 2 \psi_{k-1/2}+\psi_{k-1}}}{6}\,.
  \end{gather*}
  For the other term, on each mesh interval $I_k$ one can take the maximum of
  the values at the two end points and in the mid point of that interval.
  This gives
  \begin{gather*}
    \bignorm{\widebar{bU} - b\bar{U}}_\infty
      \approx \max_{k=1,\dots,N}\abs{\bigl(\widebar{bU} - b\bar{U}\bigr)_{k-1/2}}\,.
  \end{gather*}
  A Taylor analysis reveals that the errors introduced by these approximations
  are of order $3$ or higher.
\end{remark}

\section{Numerical results}

In this section we shall present some results of numerical experiments when
the error estimator of \autoref{theo:main} is applied to the following two test problems:
\begin{gather}\label{test}
  -\eps u'' - (2+x)u' + (1+\cos x) u = \E^{1-x}, \ \ u(0)=u(1)=0,
\end{gather}
and
\begin{gather}\label{test2}
  -\eps u'' - (2+x)u' + (1+\cos x) u = (1-x)^{\alpha}\sin{x}, \ \ u(0)=u(1)=0,
\end{gather}
where \mbox{$\alpha\in(0,1)$}.
Because of the term $(1-x)^\alpha$ on the right-hand side, the second problem
has an additional difficulty: a weak singularity at $x=1$.

The exact solution to these two problems is not available.
Therefore, we approximate the error $\chi_N$ by comparing the approximation
obtained on the mesh $\omega$ with the approximation obtained on a mesh that
is $4$ times finer than $\omega$. 

\subsection{Layer adapted meshes}

We start with examining two meshes that are \emph{a priori} designed to resolve
the layer present in the solution of our problems.

\paragraph{Shishkin meshes~\cite{Shi92}} are particularly simple because they
are piecewise uniform.
They are constructed as follows.
Let \mbox{$\sigma > 0$} be a user chosen mesh parameter.
A mesh transition point $\tau$ is defined by 
\begin{gather*}
  \tau= \min \left\{\frac{1}{2}, \frac{\sigma \eps}{\beta}\ln N \right\}.
\end{gather*}
Then the intervals \mbox{$[0,\tau]$} and \mbox{$[\tau,1]$} are each divided
into $N/2$ equidistant subintervals.

\paragraph{Bakhvalov meshes~\cite{208.19103}} are better adapted to the layers,
but less simple in construction.
A mesh generating function is defined by
\begin{gather*}
  \phi(t) \coloneqq
  \begin{cases}
     \chi(t) \coloneqq -\frac{\sigma \eps}{\rho}\ln\left(1-\frac{t}{q}\right)  & t \in [0,\tau],\\[1ex]  
     \pi(t)  \coloneqq \chi(\tau)+\chi'(\tau)(t-\tau) & \text{otherwise,}
  \end{cases}
\end{gather*}
where the scaling parameter \mbox{$q\in\left(0,1\right)$} and
\mbox{$\sigma>0$} are user chosen parameters.
The transition point $\tau$ is chosen such that \mbox{$\phi\in C^1[0,1]$}.
Now, the mesh points are given by \mbox{$x_i=\phi\left(i/N\right)$}, \mbox{$i=0,1,\dots,N$}.

In~\cite[\S3.2]{MR2046176} it was shown that the extrapolated upwind scheme is
uniformly convergent with respect to the perturbation parameter $\eps$ and that
the maximum-norm errors satisfy
\begin{gather*}
  \norm{u-U}_\infty \le
  \begin{cases}
     \left(N^{-1}\ln N\right)^2 & \text{for the Shishkin mesh and}\\[1ex]  
     N^{-2}                     & \text{for the Bakhvalov mesh,}
  \end{cases}
\end{gather*}
provided the data possesses certain regularity.

\hyperref[Results-Shishkin]{Tables \ref{Results-Shishkin}}
and \ref{Results-Bakhvalov} display the results
obtained for test problem~\eqref{test}.
The first five columns of these tables contain the number of mesh points,
the approximate error $\chi_N$, the order of convergence $p$, the value
of the a posteriori error estimator and the efficiency o the estimator,
i.e. the ratio of error and estimato.
The next five columns are reserved for the terms of the various components
of the a posteriori error estimator:
\begin{gather*}
  \eta_{\psi} \coloneqq \frac{2}{\beta}\norm{\psi-\hat{\psi}}_1,
      \quad
  \eta_{d\psi} \coloneqq \beta \max_{k=1,\dots,N}  \frac{h_k^2}{4} \abs{\psi_{x,k}},
      \quad
  \eta_{bu} \coloneqq \frac{2}{\beta}\bignorm{\widebar{bU} - b\bar{U}}_\infty,\\
  \eta_{\psi_b} \coloneqq \frac{2C^*}{\beta}  \max_{k=1,\dots,N}
            \left\{ \abs{(f-c\bar{U})_{k-1/2} + (bU)_{x,k}}\min\left[\frac{h_k}{\norm{b}_\infty},
                        \frac{h_k^2}{4\eps}\right]\right\},
       \quad
  \eta_{\Gamma\Delta} \coloneqq \frac{2}{\beta}\max_{k=1,\dots,N} \abs{\Gamma_k + \Delta_k}. 
\end{gather*} 
in order to evaluate the contribution of each part to the
a posteriori error estimator.
 
\begin{table}
  \centerline{\begin{tabular}{c|cccc|ccccc}
    $N$ & $\chi_N$ & $p$ & $\eta_N$  & $\chi_N/\eta_N$ & $\eta_{\psi}$ &  $\eta_{d\psi}$ & $\eta_{bu}$ & $\eta_{\psi_b}$ & $\eta_{\Gamma\Delta}$  \\ \hline
    $2^7$ & 1.76e-3 & 1.60 & 1.39e-2 & 1/8 & 1.51e-5 & 1.88e-5 & 6.34e-5 & 8.20e-3 & 5.59e-3  \\
    $2^8$ & 5.79e-4 & 1.65 & 4.80e-3 & 1/9 & 3.77e-6 & 4.70e-6 & 1.59e-5 & 2.77e-3 & 2.00e-3  \\
    $2^9$ & 1.84e-4 & 1.69 & 1.57e-3 & 1/9 & 9.44e-7 & 1.18e-6 & 4.00e-6 & 8.93e-4 & 6.71e-4  \\
 $2^{10}$ & 5.70e-5 & 1.72 & 4.95e-4 & 1/9 & 2.36e-7 & 2.94e-7 & 1.00e-6 & 2.79e-4 & 2.14e-4  \\
 $2^{11}$ & 1.73e-5 & 1.75 & 1.51e-4 & 1/9 & 5.90e-8 & 7.35e-8 & 2.50e-7 & 8.48e-5 & 6.61e-5  \\
 $2^{12}$ & 5.14e-6 & 1.77 & 4.53e-5 & 1/9 & 1.47e-8 & 1.84e-8 & 6.26e-8 & 2.53e-5 & 1.99e-5  \\
 $2^{13}$ & 1.51e-6 & 1.79 & 1.33e-5 & 1/9 & 3.69e-9 & 4.59e-9 & 1.57e-8 & 7.44e-6 & 5.87e-6  \\
 $2^{14}$ & 4.37e-7 & ---  & 3.88e-6 & 1/9 & 9.22e-10 & 1.15e-9 & 3.91e-9 & 2.16e-6 & 1.71e-6  \\
  \hline
  \end{tabular}}
  \caption{Test problem~\eqref{test} with $\eps=10^{-6}$,
     Shishkin mesh, $\sigma=2$.}\label{Results-Shishkin}
\end{table}

\begin{table}
  \centerline{\begin{tabular}{c|cccc|ccccc}
	$N$ & $\chi_N$ & $p$ & $\eta_N$  & $\chi_N/\eta_N$ & $\eta_{\psi}$ &  $\eta_{d\psi}$ & $\eta_{bu}$ & $\eta_{\psi_b}$ & $\eta_{\Gamma\Delta}$  \\ \hline
    $2^7$ & 1.06e-4 & 1.99 & 7.13e-4 & 1/7 & 1.51e-5 & 1.86e-5 & 6.36e-5 & 3.94e-4 & 2.21e-4  \\
    $2^8$ & 2.66e-5 & 2.00 & 1.79e-4 & 1/7 & 3.77e-6 & 4.68e-6 & 1.60e-5 & 9.82e-5 & 5.63e-5  \\
    $2^9$ & 6.65e-6 & 2.00 & 4.48e-5 & 1/7 & 9.44e-7 & 1.17e-6 & 4.00e-6 & 2.45e-5 & 1.42e-5  \\
 $2^{10}$ & 1.66e-6 & 2.00 & 1.12e-5 & 1/7 & 2.36e-7 & 2.94e-7 & 1.00e-6 & 6.11e-6 & 3.57e-6  \\
 $2^{12}$ & 4.16e-7 & 2.00 & 2.80e-6 & 1/7 & 5.90e-8 & 7.35e-8 & 2.50e-7 & 1.53e-6 & 8.94e-7  \\
 $2^{12}$ & 1.04e-7 & 2.00 & 7.01e-7 & 1/7 & 1.47e-8 & 1.84e-8 & 6.26e-8 & 3.82e-7 & 2.24e-7  \\
 $2^{13}$ & 2.61e-8 & 2.32 & 1.75e-7 & 1/7 & 3.69e-9 & 4.59e-9 & 1.57e-8 & 9.54e-8 & 5.60e-8  \\
 $2^{14}$ & 5.22e-9 & ---  & 4.38e-8 & 1/9 & 9.22e-10 & 1.15e-9 & 3.91e-9 & 2.39e-8 & 1.40e-8  \\
  \hline
  \end{tabular}}
  \caption{Test problem~\eqref{test} with $\eps=10^{-6}$,
     Bakhvalov mesh, $\sigma=2$, $q=1/2$.}\label{Results-Bakhvalov}
\end{table}

As expected, the Bakhvalov mesh performs better, i.e. give more accurate
approximations, than the Shishkin.
Comparing the two tables, we observe that the components $\eta_{\psi}$,
$\eta_{d\psi}$ and $\eta_{bu}$ are identical for the two meshes, while
there significant differences for $\eta_{\psi_b}$ and $\eta_{\Gamma\Delta}$.
The latter two dominate the error estimator closely correlate with the actual
errors.
The errors are overestimated by a factor of $9$ (Shishkin) and $7$ (Bakhvalov).
This is an acceptable overestimation.
Note that for the Bakhvalov mesh with $N=2^{14}$ mesh intervals a slight
deterioration is noticeable.
This is caused be the limitations of floating point arithmetic.

\subsection{An adaptive algorithm}
\label{ssec:adapt}

Using the a posteriori estimates of the preceding section an
adaptive algorithm can be devised.
It is based on work by de Boor~\cite{MR0431606} and uses an
equidistribution principle. 
The idea is to adaptively design a mesh for which the local
contributions to the a posteriori error estimator
\begin{align*}
  \mu_i\left(U,\omega\right)
     \coloneqq  \frac{2}{\beta}\Biggl( h_i\frac{\abs{\psi_i - 2 \psi_{i-1/2}+\psi_{i-1}}}{6} +\frac{h_i^2}{8} \abs{\psi_{x,i}} & +  \abs{\bigl(\widebar{bU} - b\bar{U}\bigr)_{i-1/2}}
        + \abs{\Gamma_i + \Delta_i} \\
       & 
          + C^*  \abs{(f-c\bar{U})_{i-1/2} + (bU)_{x,i}}
                               \min\left\{\frac{h_i}{\norm{b}_\infty},
                                         \frac{h_i^2}{4\eps}
                                   \right\}
             \Biggr),
\end{align*}
are the same on each mesh interval, i.e. \
$\mu_{i-1}\left(U,\omega\right) = \mu_i\left(U,\omega\right)$,
$i=1,\dots,N$.
This is equivalent to
\begin{gather}\label{equi}
  Q_i\left(U,\omega\right)
          = \frac{1}{N} \sum_{j=1}^{N} Q_j\left(U,\omega\right),
   \quad Q_i\left(U,\omega\right) \coloneqq
             \mu_i\left(U,\omega\right)^{1/2},
\end{gather}
However, de Boor's algorithm,
which we are going to describe now, becomes numerically unstable
when the equidistribution principle~\eqref{equi} is enforced
strongly. Instead, we shall stop the algorithm when
\begin{gather*}
  \tilde{Q}_i\left(U,\omega\right)
    \le \frac{\gamma}{N} \sum_{j=1}^{N} \tilde{Q}_j\left(U,\omega\right),
\end{gather*}
for some user chosen constant \mbox{$\gamma>1$} (in our experiments, we have
chosen \mbox{$\gamma=1.2$}).
Here we have also modified $Q_i$ by choosing
\begin{gather*}
  \tilde{Q}_i\left(U,\omega\right) \coloneqq
             \left(h_i^2+\mu_i\left(U,\omega\right)\right)^{1/2},
\end{gather*}
Adding this constant floor to $\mu_i$ avoids mesh starvation and
smoothes the convergence of the adaptive mesh algorithm.

\paragraph{\textbf{Algorithm} (de Boor,~\cite{MR0431606})}
\begin{enumerate*}
  \item Fix $N$, $r$ and a constant $\gamma>1$.
  The initial mesh $\omega^{[0]}$ is uniform with mesh size $1/N$.

  \item For $k=0, 1, \dots$, given the mesh $\omega^{[k]}$,
  compute the discrete solution $U_{\omega^{[k]}}^{[k]}$ on this mesh using
  extrapolated upwind method.
  Set $h_i^{[k]} = x_i^{[k]} - x_{i-1}^{[k]}$ for each $i$.
  Compute the piecewise-constant monitor function $M^{[k]}$ defined by
  \begin{gather*}
     M^{[k]} (x) \coloneqq \dfrac{\tilde{Q}_i^{[k]}}{h_i^{[k]}} \coloneqq \dfrac{\tilde{Q}_i\left(U_{\omega^{[k]}}^{[k]},\omega^{[k]}\right)}{h_i^{[k]}}
      \quad\text{for}\quad x\in\bigl(x^{[k]}_{i-1},x^{[k]}_i\bigr).
  \end{gather*}
  Set
  \begin{gather*}
    I^{[k]} \coloneqq \sum_{j=1}^N \tilde{Q}_j^{[k]}.
  \end{gather*}
  \item  Test the mesh:
  If
  \begin{gather*}
    \tilde{Q}_j^{[k]} \le \gamma I^{[k]} N^{-1} \quad \text{for all}
    \ \ j=1,\dots,N
  \end{gather*}
  then go to Step 5.
  Otherwise, continue to Step 4.
  \item Generate a new mesh by equidistributing the monitor
  function~$M^{[k]}$, i.e., choose the new mesh $\Delta^{[k+1]}$ such that
  \begin{gather*}
    \int_{x^{[k+1]}_{i-1}}^{x^{[k+1]}_i} M^{[k]}(t)\D t =
       \frac{I^{[k]}}N, \quad i=1,\dots,N.
  \end{gather*}
  Return to Step 2.

  \item Set $\omega = \omega^{[k]}$ and
  $U = U_{\omega^{[k]}}^{[k]}$ then stop.
\end{enumerate*}

\autoref{Results-adaptive} displays the numerical results obtained for
the adaptive algorithm when applied to our first test problem.
Additionally, the last column contains the number of iterations required by
the de Boor algorithm to meet the stopping criterion.
The accuracy achieved is comparable to that of the Bakhvalov mesh.
The rates of convergence vary because the equidistribution principle is not
imposed strictly, but they average at about $2.01$.

\begin{table}[b!]
  \centerline{\begin{tabular}{c|cccc|ccccc|c}
	$N$ & $\chi_N$ & $p$ & $\eta_N$  & $\chi_N/\eta_N$ & $\eta_{\psi}$ &  $\eta_{d\psi}$ & $\eta_{bu}$ & $\eta_{\psi_b}$ & $\eta_{\Gamma\Delta}$ & iter  \\ \hline   
      $2^7$ & 9.84e-5 & 1.76 & 7.12e-4 & 1/8 & 1.79e-5 & 2.13e-5 & 6.51e-5 & 4.53e-4 & 1.54e-4 &       6  \\
      $2^8$ & 2.90e-5 & 2.17 & 1.90e-4 & 1/7 & 4.91e-6 & 5.82e-6 & 1.80e-5 & 1.21e-4 & 3.97e-5 &       5  \\
      $2^9$ & 6.47e-6 & 2.04 & 4.71e-5 & 1/8 & 1.20e-6 & 1.41e-6 & 4.40e-6 & 2.96e-5 & 1.04e-5 &       5  \\
   $2^{10}$ & 1.57e-6 & 1.74 & 1.11e-5 & 1/8 & 2.81e-7 & 3.30e-7 & 1.03e-6 & 6.96e-6 & 2.45e-6 &       5  \\
   $2^{11}$ & 4.72e-7 & 2.17 & 2.94e-6 & 1/7 & 7.59e-8 & 8.91e-8 & 2.79e-7 & 1.87e-6 & 6.21e-7 &       4  \\
   $2^{12}$ & 1.05e-7 & 2.05 & 7.34e-7 & 1/7 & 1.89e-8 & 2.21e-8 & 6.94e-8 & 4.66e-7 & 1.58e-7 &       4  \\
   $2^{13}$ & 2.53e-8 & 2.11 & 1.85e-7 & 1/8 & 4.70e-9 & 5.52e-9 & 1.73e-8 & 1.16e-7 & 4.09e-8 &       4  \\
   $2^{14}$ & 5.85e-9 & ---  & 4.61e-8 & 1/8 & 1.17e-9 & 1.37e-9 & 4.30e-9 & 2.90e-8 & 1.03e-8 &       4  \\
  \hline
  \end{tabular}}
  \caption{Test problem~\eqref{test} with $\eps=10^{-6}$,
     adaptive algorithm.}\label{Results-adaptive}
\end{table}

Finally, we consider the second test problem, i.e.~\eqref{test2}.
As mentioned before, the term $(1-x)^\alpha$ on the right-hand side
causes a weak singularity at \mbox{$x=1$}.
Consequently, a mesh that resolves the layers only, but not the singularity,
will give inferior approximations.
This is confirmed by the results in Tables~\ref{Results1} and \ref{Results2}.
The order of convergence for the Bakhvalov mesh is reduced to
\mbox{$p=1+\alpha < 2$}, while the adaptive algorithm gives significantly
better approximation.
Nonetheless, the rates of convergence are smaller than $2$, but aproaching
$2$ as $N$ increases.

\begin{table}[h!]
\centerline{\begin{tabular}{r|cccc|cccc}
& \multicolumn{4}{|c|}{Bakhvalov mesh} &  \multicolumn{4}{|c}{adaptive algorithm}\\ 
\hline
	$N$ & $\chi_N$ & $p$ & $\eta_N$  & $\chi_N/\eta_N$ &  $\chi_N$ & $p$ & $\eta_N$  & $\chi_N/\eta_N$  \\ \hline
    $2^7$ & 1.17e-4 & 1.24 & 1.31e-3 & 1/12 &  4.92e-5 & 1.57 & 6.18e-4 & 1/13 \\
    $2^8$ & 4.94e-5 & 1.25 & 5.24e-4 & 1/11 &  1.65e-5 & 1.55 & 1.75e-4 & 1/11 \\
    $2^9$ & 2.08e-5 & 1.25 & 2.13e-4 & 1/11 &  5.63e-6 & 1.63 & 4.81e-5 & 1/9  \\
 $2^{10}$ & 8.76e-6 & 1.25 & 8.76e-5 & 1/11 &  1.81e-6 & 1.71 & 1.37e-5 & 1/8\\
 $2^{11}$ & 3.68e-6 & 1.25 & 3.63e-5 & 1/10 &  5.54e-7 & 1.74 & 3.85e-6 & 1/7 \\
 $2^{12}$ & 1.55e-6 & 1.25 & 1.51e-5 & 1/10 &  1.66e-7 & 1.82 & 1.07e-6 & 1/7\\
 $2^{13}$ & 6.50e-7 & 1.25 & 6.31e-6 & 1/10 &  4.72e-8 & 1.93 & 2.90e-7 & 1/7 \\
 $2^{14}$ & 2.72e-7 & ---  & 2.64e-6 & 1/10 &  1.24e-8 & ---  & 7.58e-8 & 1/7\\
 \hline
\end{tabular}}
  \caption{Test problem~\eqref{test2} with $\eps=10^{-6}$ and $\alpha=1/4$.}\label{Results1}
\end{table}

\begin{table}[h!]
\centerline{\begin{tabular}{r|cccc|cccc}
& \multicolumn{4}{|c|}{Bakhvalov mesh} &  \multicolumn{4}{|c}{adaptive algorithm}\\ 
\hline
	$N$ & $\chi_N$ & $p$ & $\eta_N$  & $\chi_N/\eta_N$ &  $\chi_N$ & $p$ & $\eta_N$  & $\chi_N/\eta_N$  \\ \hline
    $2^7$ & 1.04e-4 & 1.08 & 2.25e-3 & 1/22 &  4.49e-5 & 1.36 & 7.11e-4 & 1/16 \\
    $2^8$ & 4.91e-5 & 1.09 & 9.95e-4 & 1/21 &  1.75e-5 & 1.47 & 2.04e-4 & 1/12 \\
    $2^9$ & 2.30e-5 & 1.10 & 4.48e-4 & 1/20 &  6.31e-6 & 1.64 & 5.72e-5 & 1/10  \\
 $2^{10}$ & 1.08e-5 & 1.10 & 2.04e-4 & 1/19 &  2.03e-6 & 1.61 & 1.70e-5 & 1/9\\
 $2^{11}$ & 5.03e-6 & 1.10 & 9.30e-5 & 1/19 &  6.64e-7 & 1.67 & 4.84e-6 & 1/8 \\
 $2^{12}$ & 2.35e-6 & 1.10 & 4.25e-5 & 1/19 &  2.08e-7 & 1.82 & 1.35e-6 & 1/7\\
 $2^{13}$ & 1.09e-6 & 1.11 & 1.95e-5 & 1/18 &  5.90e-8 & 1.94 & 3.59e-7 & 1/7 \\
 $2^{14}$ & 5.08e-7 & ---  & 8.91e-6 & 1/18 &  1.54e-8 & ---  & 9.28e-8 & 1/7 \\\hline
\end{tabular}}
  \caption{Test problem~\eqref{test2} with $\eps=10^{-6}$ and $\alpha=1/10$.}\label{Results2}
\end{table}

\section{Conclusions}

We have studied maximum-norm \emph{a posteriori} error estimates for a
singularly perturbed convection-diffusion problem.
The differential equations has been discretised using simple upwind
difference schemes on two nested meshes.
Then these first-order approximation were combined using Richardson
extrapolation to give an enhanced method that is of second order.

The error estimators have been test for two standard layer-adapted meshes.
Furthermore, an adaptive procedure was devised that automatically detects
the location of the layer and other features of the solution (like a weak
singularity).

An interesting open question remains the derivation of \emph{a posteriori}
error bounds for higher-order extrapolation schemes that combine more than two
first-order upwind approximations.

\pagebreak

\def\cprime{$'$}

\end{document}